\documentclass[10pt]{article}

\usepackage{amssymb,latexsym, amsmath}
\usepackage{amscd}

\usepackage{amssymb}
\usepackage{amsthm}
\usepackage{amsmath}
\usepackage{amsfonts}

\newtheorem{remark}{Remark}

\newtheorem{definition}{Definition}
\newtheorem{proposition}{Proposition}
\newtheorem{theorem}{Theorem}
\newtheorem{lemma}{Lemma}

\newtheorem{example}{Example}

\newcommand\ad{\mathrm{ad}}
\newcommand\Ad{\mathrm{Ad}}
\newcommand\goth{\mathfrak}

\newcommand\ind{\mathrm{ind}\,}

\newcommand\Ker{\mathrm{Ker}\, }
\newcommand\Ann{\mathrm{Ann}\, }

\newcommand\rank{\mathrm{rank}\, }
\newcommand\corank{\mathrm{corank}\, }
\newcommand\codim{\mathrm{codim}\, }
\newcommand\Sing{{\mathsf{Sing}}}
\newcommand\Inv{Y_{\mathrm{formal}}(\goth g, a)}
\newcommand\kfield{{\mathbb{K}}}

\title{Some remarks about Mishchenko-Fomenko subalgebras}
\author{Alexey Bolsinov}

\begin{document}
\maketitle

\begin{abstract}
We discuss and compare two different approaches to the notion of Mishchenko--Fomenko subalgebras in Poisson-Lie algebras of finite-dimensional Lie algebras. One of them, commonly accepted by the algebraic community, uses polynomial $\Ad^*$-invariants. The other is based on formal $\Ad^*$-invariants and allows one to deal with arbitrary Lie algebras, not necessarily algebraic. In this sense, the latter is more universal.
\end{abstract}

\section{Motivation}\label{motivation}

This note is primarily motivated by the paper by A.~Ooms \cite{Ooms2} in which, among other interesting results, the author constructs a counterexample to my completeness criterion for Mishchenko-Fomenko subalgebras \cite{Bols1}.   I do not intend to disprove this statement by Ooms.  My point is that the example by A.~Ooms and the completeness criterion from
\cite{Bols1} are both correct.  The confusion is caused by the fact that the definitions of Mishchenko-Fomenko subalgebras used in \cite{Ooms2} and \cite{Bols1} are different.   The purpose of the present paper is to clarify this issue and perhaps to convince the reader that the definition from  \cite{Bols1} is slightly better.

\section{Formal $\Ad^*$--invariants}

To emphasise the algebraic nature of all the constructions in this paper, in what follows we consider finite-dimensional  Lie algebras over an algebraically closed field $\kfield$ of characteristic zero. However, in this section, it is not important for the field to be algebraically closed.
 
Here we recall and slightly modify the results of \cite{bz}.
I suspect that these results are not essentially new as the main issue below seems to be quite natural and was presumably previously discussed in different situations.  Basically,  we want to develop some algebraic techniques allowing us  to deal with arbitrary Lie algebras, not necessarily algebraic.  So we do not assume the existence of any polynomial and even rational $\Ad^*$-invariants.  Moreover, we never use the Lie group $G$ associated with $\goth g$.

Let $\goth g$ be a finite dimensional  Lie algebra, $\goth g^*$ its dual space and $P(\goth g)$  denote the algebra of polynomials\footnote{$P(g)$ as a set is of course the same as the symmetric Lie algebra $S(\goth g)$, but we use a slightly different point of view thinking of $P(\goth g)$ as a Poisson algebra and of its elements as functions on the vector space $\goth g^*$.} on $\goth g^*$. The algebra $P(\goth g)$ is endowed with the standard Lie-Poisson bracket

\begin{equation}
\label{LPb}
\{ f(x), g(x)\} = \langle x , [df(x), dg(x)]\rangle, \quad x\in\goth g^*,  \  df(x), dg(x)\in \goth g,
\end{equation}
and we will refer to $P(\goth g)$ as the Lie-Poisson algebra associated with $\goth g$.

Our goal is to construct a ``big'' commutative subalgebra in $P(\goth g)$.   The argument shift method suggested by A.~Mishchenko and A.~Fomenko \cite{MischFom} is based on some nice properties of $\Ad^*$-invariants.  In general, however,  polynomial (and even rational) invariants do not necessarily exist.  To avoid this problem on can use formal invariants which can be defined in the following way.

\begin{definition}
\label{def1}{\rm
Let $F= \sum_{k=1}^\infty f^{(k)}$ be a formal power series where $f^{(k)}\in P(\goth g)$ is a homogeneous polynomial of degree $k\in\mathbb N$. We say that $F$ is a {\it formal $\Ad^*$-invariant} at  a point $a\in\goth g^*$, if the following (formal) identity holds for all $\xi\in \goth g$:
\begin{equation}
\label{forminv}
\langle dF(x), \ad^*_\xi (a+x) \rangle = 0.
\end{equation}
}\end{definition}

From the differential-geometric point of view this condition simply means that the differential of $F$ at the point $a+x$ vanishes on the tangent space of the coadjoint orbit through this point.  Thus,  the above relation can be understood as the standard  definition of an invariant function $F$ where $F$ is replaced by its Taylor expansion at the point $a\in \goth g^*$.
The formal identity \eqref{forminv}  amounts to the following infinite sequence of polynomial relations:

$$
\begin{aligned}
\langle df^{(1)} (x), \ad^*_\xi a \rangle &= 0,\\
\langle df^{(2)} (x), \ad^*_\xi a \rangle &= - \langle df^{(1)} (x), \ad^*_\xi x \rangle, \\
\langle df^{(3)} (x), \ad^*_\xi a \rangle &= - \langle df^{(2)} (x), \ad^*_\xi x \rangle, \\
 \dots & \\
\langle df^{(k)} (x), \ad^*_\xi a \rangle &= - \langle df^{(k-1)} (x), \ad^*_\xi x \rangle, \\
\dots &
\end{aligned}
$$
where $\xi\in\goth g$ is arbitrary, or equivalently

\begin{equation}
\label{MFrelations}
\begin{aligned}
\ad^*_{df^{(1)} (x)} a &= 0,\\
 \ad^*_{df^{(2)} (x)} a  &= -   \ad^*_{df^{(1)} (x)} x, \\
 \ad^*_{df^{(3)} (x)} a  &= -   \ad^*_{df^{(2)} (x)} x, \\
 \dots & \\
 \ad^*_{ df^{(k)} (x)} a  &= -  \ad^*_{ df^{(k-1)} (x)} x, \\
\dots &
\end{aligned}
\end{equation}

The first relation means that the differential $df^{(1)}$ of the first term belongs to the $\ad^*$-stationary subalgebra of $a\in \goth g^*$ or equivalently:

$$
f^{(1)} \in \Ann (a) = \{ \eta\in\goth g~|~ \ad^*_\eta a =0\},
$$
as $f^{(1)}$ is a linear function and hence we may identify $f^{(1)}$ with $df^{(1)}$.

Let us denote the space of all formal $\Ad^*$-invariants at $a\in\goth g^*$ by $\Inv$. It is easy to see that this set is closed under addition and multiplication  (clearly, the usual multiplication of formal power series is well defined in our case).  Thus, any polynomial $p(F_1, \dots, F_s)$  in formal invariants $F_1,\dots, F_s\in \Inv$ is still a formal invariant.
Moreover,  if we consider a formal power series

$$
P(F_1, \dots, F_s)= \sum_{k=1}^\infty p^{(k)} (F_1, \dots, F_s)
$$
of formal invariants $F_1,\dots, F_s\in \Inv$, then $P(F_1, \dots, F_s)\in \Inv$.

The next theorem is a {\it formal} analog of some well-known facts about local invariants of a smooth action of a Lie group at a generic point.

As usual, we say that $a\in \goth g^*$ is {\it regular},  if  $\dim \Ann (a)$ is minimal and is equal to $s=\ind \goth g$.

\begin{theorem}
\label{th1}
Let $a\in\goth g^*$ be regular and $\xi_1,\dots, \xi_s$ be a basis of $\Ann (a)$.
 Then there exist formal $\Ad^*$-invariants $F_1, \dots, F_s$  such that their linear terms $f_1^{(1)}, \dots, f_s^{(1)}$ coincide with $\xi_1,\dots, \xi_s$ respectively.  Moreover, any other formal invariant $F\in\Inv$  can uniquely be represented as a formal power series of $F_1, \dots, F_s$.
\end{theorem}

\proof
The ``existence part'' of Theorem \ref{th1} was proved in \cite{bz}.   We only need to comment on the second part.  This statement immediately follow from

\begin{lemma}
\label{lem1}
 Suppose that $F\in \Inv$ starts with a term of degree $m$, i.e., $F = \sum_{k=m}^\infty f^{(k)}$.  Then there is a homogeneous polynomial $p^{(m)}(F_1, \dots, F_s)$ of degree $m$  such that $F - p^{(m)}(F_1, \dots, F_s)$ starts with a term of degree $m+1$.
\end{lemma}

\proof
Since $F=\sum_{k=m}^\infty f^{(k)}$ is a formal invariant, we have a sequence of relations similar to \eqref{MFrelations} but these relations start with the identity
$$
\langle df^{(m)} (x), \ad^*_\xi a \rangle = 0  \quad \mbox{for all } \xi\in\goth g,
$$
or equivalently,
$$
\ad^*_{df^{(m)} (x)} a  = 0,
$$
i.e.,  $df^{(m)} (x) \in \Ann(a)$.  In other words, $f^{(m)}$ is a homogeneous polynomial of degree $m$ on $\goth g^*$ whose differential at every point $x\in \goth g^*$ belongs to
$\Ann(a)$.  But this condition obviously means that $f^{(m)}$ can, in fact, be written as a polynomial of the basis elements $\xi_1, \dots, \xi_s\in \Ann (a)$, i.e. $f^{(m)} = p^{(m)}(\xi_1, \dots, \xi_s)$.

Since $\xi_i$'s are the linear terms of $F_i$'s,   we see that the $m$-terms of  $F$ and  $p^{(m)}(F_1, \dots, F_m)$ coincide so that the power series $F- p^{(m)}(F_1, \dots, F_m)$ starts from a term of degree $m+1$, as required. \endproof

Thus, Lemma \ref{lem1} says the following. Given a formal invariant $F=\sum f^{(k)} \in\Inv$ we can step by step kill all of its homogeneous terms by subtracting a suitable polynomial  $p^{(m)}(F_1, \dots, F_s)$, $m=1,2,\dots$, in other words, $F = \sum p^{(m)}(F_1, \dots, F_s)$ as required.

The uniqueness of such an expansion follows from the independence of $F_1, \dots, F_s$ so that no nontrivial polynomial in $F_1, \dots, F_s$ may vanish identically.  \endproof

\begin{remark}
{\rm A similar result, of course, holds true for an arbitrary linear representation $\rho : \goth g \to \mathrm{End}(V)$  (see \cite{bz}).  The only difference is that $\Ann (a)$ should be replaced by the orthogonal complement  $T_a^\bot\subset V^*$ to the ``$\rho$-orbit''     $T_a = \{ \rho (\xi)  a, \ \xi\in\goth g\}\subset V$. The natural identification $\Ann (a) \simeq T_a^\bot$ makes no sense and does not have any analog  in the case when $\rho \not\simeq \ad^*$.

}\end{remark}

Also, it is interesting to notice that finding a formal invariant  (up to any order) is a problem of Linear Algebra.  Indeed,  we only need to solve successively the sequence of relations \eqref{MFrelations}.  The ``existence  part'' of Theorem \ref{th1}  tells us that the systems of linear equations we obtain at each step are all consistent and, moreover, the proof given in \cite{bz} explains how to make the choice of a solution $f^{(k)}$ unique.

\section{Mishchenko--Fomenko subalgebras: two versions}

We first recall the definition used in \cite{bolsActa, Bols1}.

\begin{definition} 
\label{MF2}{\rm
Let $a\in \goth g^*$ be regular  and $F_1=\sum f_1^{(k)}, \dots, F_s=\sum f_s^{(k)} \in \Inv$ be an arbitrary basis of formal $\Ad^*$-invariants at the point $a\in \goth g^*$ as in Theorem \ref{th1}.
The {\it algebra of polynomial shifts} $\mathcal F_a(\goth g)$ is defined to be the subalgebra in $P(\goth g)$ generated by the homogeneous polynomials $f_i^{(k)}$, $i=1,\dots, s=\ind \goth g$, $k\in\mathbb N$.
}\end{definition}

\begin{remark}{\rm
In \cite{bolsActa, Bols1},  instead of formal $\Ad^*$-invariants we considered the Taylor expansions of smooth or  analytic (local) invariants $F$ which always exist in a neighbourhood of a regular point $a\in \goth g^*$:
$$
F(a + t x) \simeq \sum_{k=0}^\infty t^k f^{(k)} (x).
$$
The above definition is just a straightforward extension (or algebraic reformulation) of this construction to the case of an arbitrary field of characteristic zero.   The term ``algebra of polynomial shifts'' was used in our recent paper \cite{BolsZhang}  to emphasise the difference from ``standard shifts''  $F(x+ t a)$ which are not necessarily polynomial in $x$ (e.g., if $F$ itself is not a polynomial).   We consider $\mathcal F_a(\goth g)$ as the first version of a Mishchenko--Fomenko subalgebra but do not use this terminology for $\mathcal F_a(\goth g)$ to avoid further confusion with another type of Mishchenko--Fomenko subalgebras discussed below.
}\end{remark}

The following proposition summarises the main properties of the algebra of polynomial shifts.

\begin{theorem}\label{th2}
Let $a\in \goth g ^*$ be an arbitrary regular element and $\mathcal F_a(\goth g)$ the corresponding algebra of polynomial shifts. Then the following properties hold:
\begin{enumerate}

\item $\mathcal F_a(\goth g)$ does not depend on the choice of the basis formal invariants $F_1, \dots, F_s\in\Inv$.

\item The linear polynomials from $\mathcal F_a(\goth g)$ are elements of $\Ann (a)$. In other words, $\mathcal F_a(\goth g) \cap \goth g = \Ann(a)$.

\item $\mathcal F_a(\goth g)$ is commutative  w.r.t. to the standard Lie-Poisson bracket \eqref{LPb}.

\item $\mathcal F_a(\goth g)$ is commutative  w.r.t. to the Poisson $a$-bracket
$$
\{ f, g\}_a =  \langle a , [df(x), dg(x)]\rangle, \quad x\in\goth g^*,  \  df(x), dg(x)\in \goth g.
$$

\item $\mathcal F_a(\goth g)$  is complete, i.e.,  $\mathrm{tr.deg.}\, \mathcal F_a(\goth g) = \frac{1}{2} (\dim \goth g + \ind\goth g)$, if and only if $\codim \Sing \ge 2$, where
$$
\Sing=\{ y \in\goth g^*~|~  \dim\Ann y >\ind \goth g\} \subset \goth g^*
$$
is the set of singular points in $\goth g^*$.

\item In general, the number of algebraically independent polynomials in $\mathcal F_{a}$ is
$$
\mathrm{tr.deg.}\, \mathcal F_a(\goth g) = \frac{1}{2} (\dim \goth g + \ind\goth g) - \mathrm{deg}\, \mathsf{p}_{\goth g},
$$
where $\mathsf{p}_{\goth g}$ is the fundamental semi-invariant of $\goth g$.
\end{enumerate}

\end{theorem}

\begin{remark}
{\rm For item 5,  it is important that $\kfield$ is algebraically closed. The other items hold true without this assumption.}
\end{remark}

\proof
Items 2, 3, 4 and 5  have been discussed in many papers (see, for instance, \cite{bolsActa, Bols1, BolsZhang, BO-rcd, bz}).  The item 5 is the completeness criterion from \cite{Bols1}.  The counterexample from \cite{Ooms2} uses a different definition of a Mishchenko-Fomenko algebra and does not contradict to item 5 (see Example \ref{ex1} below).  The item 6  was recently proved in \cite{JS}  for the Mishchenko--Fomenko algebras $Y_a(\goth g)$ in the sense of Definition \ref{MF2} below, but it is still true for $\mathcal F_a(\goth g)$ if $a$ is regular.  Theorem \ref{th3}  proved below immediately implies both  5 and 6.

We only need to explain item 1 which is fairly easy.   Let $F=\sum f^{(m)}\in\Inv$ be an arbitrary formal $\Ad^*$-invariant. It is sufficient to show that each term $f^{(m)}$ of this formal series belongs to the Mishchenko-Fomenko subalgebra $\mathcal F_a(\goth g)$. We know from Theorem \ref{th1} that $F$ can be written as a formal power series in $F_1, \dots, F_s$.  But this immediately implies that every term $f^{(m)}$ admits a polynomial representation via $f_i^{(k)}$, $i=1,\dots, s=\ind  \goth g$, $k\le m$,  and hence belongs to $\mathcal F_a(\goth g)$ as required.
\endproof

 The next definition of Mishchenko-Fomenko subalgebra is used in \cite{ JS, Ooms2,   PY} and seems to be more common in algebraic literature.

\begin{definition}\label{MF2}{\rm
Let  $Y(\goth g)= P(\goth g)^{\goth g} \subset P(\goth g)$ be the algebra of $\Ad^*$-invariant polynomials or, equivalently, the centre of $P(\goth g)$.  For  $f\in Y(\goth g)$, $a\in\goth g^*$ and $ t\in \kfield$ consider the expansion
\begin{equation}\label{shift}
f ( x + t a) =  \sum f_{a,m} (x) t^m
\end{equation}
into the powers of $t$. The polynomials  $f_{a,m} (x)$ are called the $a$-shifts of $f$.  The {\it Mishchenko--Fomenko algebra} $Y_a(\goth g)$ is defined as the subalgebra  in $P(\goth g)$  generated by the $a$-shifts  $f_{a,k}$  of  all $f\in Y(\goth g)$ (or equivalently of the generators of $Y(\goth g)$).
}\end{definition}

We first notice that if $a\in\goth g^*$ is regular, then  $Y_a(\goth g)$ is a subalgebra of $\mathcal F_a(\goth g)$.
Indeed, without loss of generality we may assume that the generators $f\in Y(\goth g)$  are homogeneous, then $f(x+ ta) = t^d f(a+ t^{-1} x)$, where $d=\mathrm{deg}\, f$ and the expansion \eqref{shift} is, in fact, equivalent to the Taylor expansion of $f(x)$ at the point $a$:
$$
f( a + x ) =\sum_{m=0}^d f^{(m)}, \quad \mbox{where } f^{(m)} = f_{a,d-m}.
$$

Since $\sum_{m=0}^d f^{(m)}$ is a formal invariant at the point $a\in\goth g^*$ in the sense of Definition \ref{def1},   all the $a$-shifts $f_{a,k}$ belong to $\mathcal F_a(\goth g)$ and consequently $Y_a(\goth g)\subset \mathcal F_a(\goth g)$.

On the other hand, the algebra of polynomial $\Ad^*$-invariants $Y(\goth g)$ might be trivial even if $\ind g = s > 0$.  In such a case, the Mischenko-Fomenko algebra $Y_a(\goth g)$ is trivial too in contrast to $\mathcal F_a(\goth g)$ that remains non-trivial since independent formal $\Ad^*$-invariants $F_1, \dots, F_s$ always exist (Theorem \ref{th1}).

To illustrate this phenomenon and to show how to describe $\mathcal F_a(\goth g)$ in practice, we consider  {\it Counterexample to Bolsinov's assertion} from \cite{Ooms2}.

\begin{example}\label{ex1}
{\rm
Consider  the solvable Lie algebra $\goth g$ of dimension 8 and index 2 defined by the relations:
$$
\begin{aligned}
&[x_0, x_1]=5x_1, \  [x_0, x_2]=10x_2, \ [x_0, x_3]=-13 x_3, \ [x_0, x_4]=-8x_4, \ [x_0, x_5]=-3x_5, \\
&[x_0, x_6]=2x_6, \ [x_0, x_7]=7x_7, \  [x_1,x_3]=x_4,\   [x_1,x_4]=x_5,\  [x_1,x_5]=x_6,\  [x_1,x_6]=x_7,\\
&[x_2,x_3]=x_5,\  [x_2,x_4]=x_6,\  [x_2,x_5]=x_7.
\end{aligned}
$$
The algebra of polynomial $\Ad^*$-invariants is trivial, i.e., $Y(\goth g)=\{ \kfield \}$ and therefore the Mishchenko--Fomenko subalgebra $Y_a(\goth g)$ is trivial too.  The singular set $\Sing$ has codimension 3 and is defined by three linear equations $\{ x_5=x_6=x_7=0\}$.   Thus, according to the completeness criterion from \cite{Bols1}  ({\it Bolsinov's assertion}),  the algebra of polynomial shifts $\mathcal F_a(\goth g)$,  $a\notin \Sing$,  is complete, i.e., $\mathrm{tr.deg.}\mathcal F_a(\goth g) = \frac{1}{2} (\dim \goth g + \ind \goth g) = 5$.

The Lie algebra $\goth g$ possesses two independent rational $\Ad^*$-invariants. If they are given explicitly, then the coefficients of their Taylor expansions at the point $a\in \goth g$ can be taken as generators of $\mathcal F_a(\goth g)$.  However, even if we do not have any information about them (I did not have it),  we can still use formal $\Ad^*$-invariants to construct $5$ algebraically independent polynomial shifts.  As an example, take $a\in \goth g^*$ such that $x_7(a)=1$ and $x_i(a)=0$, $i=0,\dots, 6$.
The stationary subalgebra $\Ann (a)$ is  generated by $x_3$ and $x_4$ and therefore according to Theorem \ref{th1} there exist formal invariants of the form:
$$
\begin{aligned}
F &= x_3 + f^{(2)} + f^{(3)} + \dots \\
H &= x_4 + h^{(2)} + h^{(3)} + \dots
\end{aligned}
$$
The ``higher'' terms can easily be found  successively by solving relations \eqref{MFrelations}.   Moreover, the solution is unique if in addition we require that
$f^{(i)}$ and $h^{(i)}$ vanish identically on the two-dimensional subspace defined by $x_0=x_1=x_2=x_5=x_6=x_7=0$.  Here is the result of the computation I have done by hand:
$$
\begin{aligned}
f^{(2)}=-\frac{13}{7}x_3x_7 + x_4x_6 + \frac{1}{2} x_5^2      , &\quad f^{(3)}= \frac{39}{49}  x_3 x_7^2 - \frac{3}{7} x_5^2 x_7 - \frac{6}{7} x_4x_6x_7\\
h^{(2)}= -\frac{8}{7} x_4x_7 + x_5x_6   , \quad &\quad h^{(3)}=\frac{4}{49} x_4 x_7^2 - \frac{1}{7} x_5x_6x_7 + \frac{1}{3} x_6^3.
\end{aligned}
$$

There is no need to continue this process,  as we have already found 5 algebraically independent polynomial shifts: $x_3$, $x_4$, $f^{(2)}$, $h^{(2)}$ and one of $f^{(3)}$, $h^{(3)}$.
Thus, the algebra $\mathcal F_a(\goth g)$ so obtained is complete despite the fact that $Y(\goth g)$ is trivial.

}\end{example}

This phenomenon was well understood long ago and, in fact,  was the main reason for us to slightly modify the original construction by A.Mishchenko and A.Fomenko in order to avoid the problem with non-existence of polynomial invariants and construct a commutative subalgebra of $P(\goth g)$  as large as possible.
To the best of my knowledge this modification is due to Andrey Brailov who explained this construction to me in 1986 when I was a PhD student. I am not sure, however, if he ever published this important remark.

The following proposition gives an obvious necessary and sufficient condition for $\mathcal F_a(\goth g)$ and $Y_a(\goth g)$ to coincide.

\begin{proposition}\label{prop1} 
The following conditions are equivalent:
\begin{enumerate}
\item  $\Ann (a)$ is generated by the differentials $df(a)$, $f\in Y(\goth g)$,
\item  $a\in \goth g^*$ is regular and $\mathcal F_a(\goth g) = Y_a(\goth g)$.
\end{enumerate}
\end{proposition}

\proof Notice that Condition 1 implies that $a$ is regular. Furthermore,  if  homogeneous invariants $f_1, \dots, f_s\in Y(\goth g)$ are such that $df_1(a), \dots, df_s(a)$ form a basis of $\Ann(a)$, then we may consider the Taylor expansions of $f_1, \dots, f_s$ at $a\in\goth g$
$$
f_i(a+x)= \sum_m f_i^{(m)}, \quad  i=1,\dots, s=\ind\goth g,
$$
as a basis in $\Inv$.  Since the homogeneous terms  $f_i^{(m)}$  in these expansions are the same as the $a$-shifts of $f_i$, we immediately conclude that $\mathcal F_a (\goth g)\subset Y_a(\goth g)$ and hence, $\mathcal F_a (\goth g)= Y_a(\goth g)$.

On the other hand,  assume that $a\in \goth g^*$ is regular and $\mathcal F_a(\goth g) = Y_a(\goth g)$.   Let us compare the linear functions contained in $\mathcal F_a(\goth g)$ and  $Y_a(\goth g)$.  According to item 2 of Theorem \ref{th2},   the linear functions of  $\mathcal F_a(\goth g)$ are exactly the elements of $\Ann (a)$. On the other hand, the linear functions from $Y_a(\goth g)$ are the differentials $df(a)$, $f\in Y(\goth g)$.   Since $\mathcal F_a(\goth g) = Y_a(\goth g)$, we get the desired conclusion.  \endproof

There are many examples of $\goth g$ and $a\in\goth g^*$ for which the above condition is fulfilled. The most important of them are semisimlpe (reductive) Lie algebras.

If $\mathrm{tr.deg.} Y(\goth g) < \ind \goth g$, then $Y_a (\goth g)$ is strictly smaller than $\mathcal F_a(\goth g)$.   On the contrary, if $\mathrm{tr.deg.} Y(\goth g) = \ind \goth g$, then $Y_a (\goth g)$  and $\mathcal F_a(\goth g)$ coincide for almost all regular $a\in\goth g$.  However,  if  $a$ is regular but the differentials of the polynomial invariants do not generate $\Ann (a)$, then we have proper inclusion  $Y_a(\goth g) \subsetneq \mathcal F_a(\goth g)$ (although $\mathrm{tr.deg.} Y_a(\goth g) = \mathrm{tr.deg.} \mathcal F_a(\goth g)$).

\begin{example}{\rm
Consider,  for instance,  the six-dimensional nilpotent Lie algebra $\goth g$ with relations (this is the Lie algebra $\goth g_{7,1.1(i_\lambda),\lambda = 1}$ with number 155 from the list  presented in \cite{Ooms2}\footnote{There is nothing special in this example. My choice was more or less random within a sub-list of Lie algebras with some suitable properties.}):
$$
\begin{aligned}
& [x_1,x_2]=x_3, \ [x_1,x_3]=x_4, \ [x_1,x_4]=x_5,\ [x_1,x_5]=x_6, \ [x_1,x_6]=x_7,  \\
& [x_2,x_3]=x_5, \ [x_2,x_4]=x_6, \ [x_3,x_4]=x_7.
\end{aligned}
$$

It is straightforward to verify that $\ind \goth g = 3$ and  the singular set $\Sing \subset \goth g^*$ is defined by three equations $\{ x_5=x_6=x_7=0\}$ so that $\codim\Sing=3$.
The algebra $Y(\goth g)$ of polynomial invariants is generated by four polynomials (see \cite{Ooms2}):
$$
x_7, \quad f=x_6^2 - 2 x_5x_7, \quad  g= 2 x_6^5 -10 x_5 x_6^3 x_7 + 15x_5^2x_6x_7^2 -15x_4x_5x_7^3 +15x_3x_6x_7^3 -15x_2x_7^4,
$$
and
$$
\begin{aligned}
& h=(4f^5 - g^2)/x_7^3 = -225x_2^2x_7^5+450x_2x_3x_6x_7^4-450x_2x_4x_5x_7^4+450x_2x_5^2x_6x_7^3-300x_2x_5x_6^3x_7^2+ \\
& 60x_2x_6^5x_7-225x_3^2x_6^2x_7^3+450x_3x_4x_5x_6x_7^3-450x_3x_5^2x_6^2x_7^2+300x_3x_5x_6^4x_7-60x_3x_6^6-\\
& 225x_4^2x_5^2x_7^3+450x_4x_5^3x_6x_7^2-300x_4x_5^2x_6^3x_7+60x_4x_5x_6^5-128x_5^5x_7^2+95x_5^4x_6^2x_7-20x_5^3x_6^4,
\end{aligned}
$$
which satisfy one relation $4f^5 - g^2 - h x_7^3 = 0$.

Since  $\mathrm{tr.deg.} Y(\goth g)=\ind\goth g$,  the differentials $df(a)$, $f\in Y(\goth g)$ generate $\Ann(a)$ for almost all regular points $a\in\goth g^*$ but not for all in this case.
From the point of view of the theory of integrable Hamiltonian systems, it is natural to think of the generators $x_7, f, g$ and $h$ as first integrals of a Hamiltonian system on $\goth g^*$,  and consider the momentum mapping $\Phi=(x_7, f, g, h): \goth g^* \to \kfield^4$.  At a generic point the differential of this map has rank 3 and it makes sense to introduce the set of critical points of $\Phi$
$$
\mathsf{Crit} = \bigl\{ y\in \goth g^*~|~ \dim \mathrm{span}\{ df(y), \  f\in Y(\goth g) \} < \ind \goth g\bigr\}.
$$

In the notation from \cite{JS}, the complement to this set $\mathsf{Crit}$ can be written as $\goth g^*_{\mathrm{Reg}}$ in  contrast to $\goth g^*_{\mathrm{reg}} = \goth g^* \setminus \Sing$. In our example, $\mathsf{Crit}$ is defined by two equations $\{ x_6 = x_7=0\}$ so that  $\mathsf{Crit}$ is larger than $\Sing$  (or equivalently, $\goth g^*_ {\mathrm{Reg}}$ is smaller than $\goth g^*_{\mathrm{reg}}$).

This means that there are regular elements $a\in\goth g^*$ for which $Y_a(\goth g) \subsetneq \mathcal F_a(\goth g)$, namely it is so for every $a\in\goth g^*$ with $x_6(a)=x_7(a)=0$, $x_5(a)\ne 0$.
For instance, if we take $a\in \goth g^*_{\mathrm{reg}} \setminus \goth g^*_{\mathrm{Reg}} = \mathsf{Crit}\setminus\Sing$ such that $x_5(a)=1$, $x_i(a)=0$, $i\ne 5$,  then it is easy to verify that $Y_a(\goth g)$ contains only one linear function, namely, $x_7$,  whereas $\mathcal F_a(\goth g)$ contains three: $x_5$, $x_6$ and $x_7$  (a basis of $\Ann (a)$).
}\end{example}

\section{One important property of the algebra of polynomial \\  shifts $\mathcal F_a(\goth g)$}

One of the advantages of the algebras $\mathcal F_a(\goth g)$ is a natural description of the subspace in $\goth g$ spanned by the differential of $f\in\mathcal F_a(g)$ at any point $x\in\goth g^*$. We denote this subspace by
$$
d\mathcal F_a (x) = \mathrm{span}\bigl\{ df (x), \  f\in\mathcal F_a(\goth g) \bigr\} \subset \goth g.
$$
This description is very simple and can be given in terms of the pencil of skew-symmetric forms generated by the forms
$$
\mathcal A_x (\xi, \eta) = \langle x, [\xi, \eta] \rangle  \quad \mbox{and} \quad \mathcal A_a (\xi, \eta) = \langle a, [\xi, \eta] \rangle.
$$

The following statement is well known \cite{Bols1, BO-rcd, bz}.

\begin{proposition}
\label{prop2}
$d\mathcal F_a (x)  =  \sum \Ker \mathcal A_{x+\lambda a}$,
where the sum is taken over all\footnote{It is sufficient to consider finitely many values of $\lambda$ in this sum. For example,  one can arbitrarily choose distinct rational numbers $\lambda_1, \dots, \lambda_k$ with $k=\dim \goth g$.} $\lambda \in \kfield$   such that $x+\lambda a \notin\Sing$.
\end{proposition}

Recall that a pair of skew-symmetric forms can simultaneously be reduced to an elegant Jordan-Kronecker canonical form \cite{IKozlov, Thompson} playing an important role in the theory of compatible Poisson brackets \cite{ BolsIzo, BO-rcd, gelzak, Zak1}. Here we formulate one straightforward and simple corollary of the Jordan-Kronecker decomposition theorem referring to \cite{BolsZhang} for details.

Let $A$ and $B$ be two skew-symmetric forms on a finite-dimensional vector space $V$, we will think of them as just two skew-symmetric matrices.  Let $r= \max_{\lambda \in \kfield} \rank (A+\lambda B)$ be the rank of the pencil of skew-symmetric forms $\mathcal P = \{A + \lambda B\}$.  Without loss of generality we assume that $B$ is regular in this pencil, i.e. $\rank B = r$.

Consider the Pfaffians of all $r\times r$ diagonal minors of $A+\lambda B$ as polynomials in $\lambda$ and  denote by $\mathsf p$ their greatest common divisor.   Notice that $\mathsf p =1$ if and only if the rank of $A + \lambda B$ never drops, i.e., equals $r$ for each $\lambda\in\kfield$.
The following formula is a corollary of the Jordan--Kronecker decomposition theorem.

\begin{proposition} 
Let $L= \sum \Ker (A+\lambda B)$ where the sum is taken over all $\lambda\in\kfield$ such that $\rank (A+\lambda B) = r$. Then
$
\dim L = \frac{1}{2} (\dim V + \corank \mathcal P) - \deg \mathsf p.
$
\end{proposition}

Let us transfer and apply this formula to our pencil of skew-symmetric forms $\mathcal P= \{\mathcal A_{x+\lambda a}\}$ on~$\goth g$.  Obviously,
$\dim V = \dim \goth g$ and $\corank \mathcal P = \ind \goth g$ (here we use the fact that $a\in\goth g^*$ is regular).   So  we only need to clarify the meaning of $\mathsf p$.  This (kind of a) polynomial is known as the fundamental semi-invariant $\mathsf p_{\goth g}$ of $\goth g$.   To define  $\mathsf p_{\goth g}$ consider the Pfaffians $p_1, \dots, p_N$ of all $r\times r$ diagonal minors of the matrix $\mathcal A_y = \bigl( c_{ij}^k y_k \bigr)$, $r=\dim\goth g - \ind\goth g$.  Then $\mathsf p_{\goth g}$ is the greatest common divisor of $p_1, \dots, p_N$  (all these polynomials are now considered as elements of $P(\goth g)$, i.e., as polynomials in $y_1\dots, y_n$).
Thus we have,
$$
\begin{aligned}
p_1(y) &= \mathsf{p}_{\goth g} (y) \cdot h_1(y) \\
 & \dots \\
p_N(y) &= \mathsf{p}_{\goth g} (y) \cdot h_N(y) \\
 \end{aligned}
$$
where $h_1(y),\dots, h_N(y)$ do not have any non-constant common factors.  This implies, by the way, that the singular set  $\Sing$ is the union of two subsets
$$
\Sing_0 = \{ \mathsf p_{\goth g}=0\} \quad \mbox{and} \quad \Sing_1 = \{ \mathsf h_1(y)=0, \dots, h_N(y)=0\}.
$$
Thus, there are three possibilities:
\begin{itemize}
\item $\mathsf{p}_{\goth g}= 1$ and then  $\Sing_0 =\emptyset$, $\Sing = \Sing_1$ and $\codim\Sing \ge 2$,
\item $h_i = \mathrm{const} \in \kfield$ and then $\Sing = \Sing_0$, $\codim\Sing = 1$ and $\Sing_1 = \emptyset$,
\item both $\mathsf{p}_{\goth g}$ and $h_i$ are non-constant,  then both $\Sing_0$ and $\Sing_1$ are non-empty and $\codim \Sing_0 = 1$ and $\codim\Sing_1 \ge 2$.
\end{itemize}

Replacing $y$  by  $x+\lambda a$  we obtain two possibilities:   either $\mathsf p_{\goth g} (x+\lambda a)$ is still a greatest common divisor of $p_1(x+\lambda a),\dots, p_N(x+\lambda a)$  (now we consider them as polynomials in one single variable $\lambda$), or the greatest common divisor $\mathsf p_{x,a}(\lambda)$ is ``bigger''.   The latter condition simply mean that $h_1(x+\lambda a), \dots, h_N(x+\lambda a)$ have a non-trivial common factor,  or  in geometric terms, that the straight line $x+\lambda a$, $\lambda\in\kfield$, intersects the set $\Sing_1$.

Thus, we come to the following conclusion which is similar to the Joseph--Shafrir formula (Section 7.2 in \cite{JS}).
Notice that this is a straightforward corollary of the Jordan--Kronecker decomposition theorem.

\begin{theorem}
\label{th3}
Let $a \in \goth g^*$ be regular and $d\mathcal F_a (x) = \mathrm{span}\bigl\{ df (x), \  f\in\mathcal F_a(\goth g)\bigr\} \subset \goth g$, $x\in \goth g^*$. Then
$$
\dim d\mathcal F_a (x)  = \frac{1}{2} (\dim\goth g + \ind\goth g) - \deg \mathsf{p}_{x,a},
$$
where $\mathsf{p}_{x,a}(\lambda)$ is the greatest common divisor of the
Pfaffians $p_1(x+\lambda a), \dots, p_N(x+\lambda a)$ of all $r\times r$ diagonal minors of the matrix $\mathcal A_{x+\lambda a} = \bigl( c_{ij}^k (x_k +\lambda a_k) \bigr)$, $r=\dim \goth g-\ind\goth g$.

In particular,
$$
\dim d\mathcal F_a (x)  \le \frac{1}{2} (\dim\goth g + \ind\goth g) - \deg \mathsf{p}_{\goth g},
$$
with equality if and only if the straight line $x+\lambda a$, $\lambda\in\kfield$,  does not intersect the subset $\Sing_1\subset \Sing$.
\end{theorem}

 A similar formula holds true for Mishchenko--Fomenko subalgebras $Y_a(\goth g)$ after some additional amendments.  For each $x\in\goth g^*$ consider the subspace
$d Y_a (x)=\mathrm{span}\bigl\{ df (x), \  f\in Y_a(\goth g)\bigr\} \subset \goth g$. As noticed above, for regular  $a\in \goth g^*$ we have the inclusion $Y_a(\goth g) \subset \mathcal F_a(\goth g)$ and therefore
$dY_a(x) \subset d\mathcal F_a(x)$ for any $x\in\goth g^*$.
A sufficient condition for these two subspaces to coincide is very simple (cf. Proposition \ref{prop1}).

\begin{proposition}
Let  $a\in \goth g^*$ be regular. If
the straight line $x + \lambda a$ do not belong to $\mathsf{Crit}$, then  $d Y_a (x) = d\mathcal F_a(x)$.
 \end{proposition}

  \proof 
  Indeed, if $y=x+\lambda a\notin\mathsf{Crit}\cup\Sing$, then the differentials of the shifted invariants $f_{\lambda}(x)=f(x+\lambda a)\in Y_a(\goth g)$, $f\in Y(\goth g)$, generate $\Ker \mathcal A_{x+\lambda a}$. Thus,
$\Ker \mathcal A_{x+\lambda a} \subset dY_a(x)$ for infinitely many $\lambda$'s  and in view of Proposition \ref{prop2}, we have the converse inclusion $d\mathcal F_a(x)\subset dY_a(x)$. 
\endproof

On the other hand, $Y_a(\goth g)$ is well defined for any $a\in\goth g^*$ both regular and singular, whereas $\mathcal F_a(\goth g)$ in general makes no sense for singular $a\in\goth g^*$.  Nevertheless, the description of the subspace $dY_a(x)$ is easy to obtain if we notice that
$d Y_a (x) = dY_x (a)$ and more generally this subspace $dY_x (a)$ depends only of the two-dimensional subspace of $\goth g^*$ generated by $a$ and $x$ so that $d Y_a (x) = dY_{a'} (x')$ if $\mathrm{span}(a', x')=\mathrm{span}(a, x)$.   In particular, if the straight line $x+\lambda a$ do not belong to the singular set, we may assume without loss of generality that $x$ is regular.  Then we have

\begin{proposition}\label{prop5}
Let $x\in\goth g^*$ be regular and the straight line 
$a + \lambda x$ do not belong to $\mathsf{Crit}$. Then $dY_a (x) =  d\mathcal F_x (a)$.
\end{proposition}

Hence we immediately obtain the following version of Theorem \ref{th3} for the Mishchenko--Fomenko subalgebras $Y_a(\goth g)$ (simply by interchanging $x$ and $a$).

\begin{theorem}
\label{th4}
Let $x\in \goth g^*$ be regular and $d Y_a (x) = \mathrm{span}\bigl\{df (x), \  f\in Y_a(\goth g)\bigr\} \subset \goth g$.  Assume that $\mathrm{tr.deg.} Y(\goth g)=\ind \goth g$ and the straight line $a+\lambda x$  does not belong to $\mathsf{Crit}$.  Then
$$
\dim dY_a (x)  = \frac{1}{2} (\dim\goth g + \ind\goth g) - \deg \mathsf{p}_{a,x},
$$
where $\mathsf{p}_{a,x}(\lambda)$ is the greatest common divisor of the
Pfaffians $p_1(a+\lambda x), \dots, p_N(a+\lambda x)$ of all $r\times r$ diagonal minors of the matrix $\mathcal A_{a+\lambda x} = \bigl( c_{ij}^k (a_k +\lambda x_k) \bigr)$, $r=\dim \goth g-\ind\goth g$.

Furthermore,
$$
\dim d Y_a (x)  \le \frac{1}{2} (\dim\goth g + \ind\goth g) - \deg \mathsf{p}_{\goth g},
$$
with equality if and only if the straight line $a+\lambda x$, $\lambda\in\kfield$,  does not intersect the subset $\Sing_1\subset \Sing$.

For a fixed $a\in \goth g^*$ such a line exists if and only if  $a\notin \Sing_1$. In particular,
$$
\mathrm{tr.deg.} Y_a (\goth g)  \le \frac{1}{2} (\dim\goth g + \ind\goth g) - \deg \mathsf{p}_{\goth g},
$$
with equality if and only if $a\notin\Sing_1$.

\end{theorem}

\begin{remark}
{\rm
The latter statement of this theorem is the Joseph--Shafrir formula (Section 7.2 in \cite{JS}). In particular,  $\Sing_1$ must coincide with the set $\goth g^* \setminus \goth g^*_{\mathrm{wreg}}$ from \cite{JS}. Our definition of  $\Sing_1$ seems to be simpler and more transparent than that of $\goth g^*_{\mathrm{wreg}}$ in \cite{JS}. Unfortunately, I was not able to verify the equivalence of these two definitions.
}\end{remark}

\end{document}